\documentclass[12pt]{article}
\usepackage[T2A]{fontenc}
\usepackage[cp1251]{inputenc}
\usepackage{amssymb,amsmath,latexsym,amsthm}
\usepackage{indentfirst}
\usepackage{verbatim}
\theoremstyle{plain}

\newtheorem{theorem}{Теорема}

\theoremstyle{remark}

\DeclareMathOperator{\Hol}{Hol}

\title{
{\large \bf Holomorphic minorants\\ of (pluri-)subharmonic functions\footnote{Supported by RFBR (project no.~13-01-00030)}}}
\author{Timur Baiguskarov, Bulat Khabibullin}
\begin{document}
\maketitle          

\section{{General case}} As usual, $\mathbb N:=\{1,2, \dots\}$, $\mathbb R$ and $\mathbb C$ are the sets of real and complex numbers, resp.
We denote by $\lambda$  the Lebesgue measure on $\mathbb C^n$, $n\in \mathbb N$. Given  $z\in \mathbb C^n$ and $r>0$, 
$B(z,r):=\bigl\{z'\in \mathbb C^n \colon |z'-z|<r\bigr\}$ . For a function $f \colon B(z,r)\rightarrow [-\infty,+\infty]:=\{-\infty\}\cup \mathbb R \cup \{+\infty\}$, we set
\begin{equation*}
	B_f (z,r):=\frac{1}{\lambda \bigl(B(z,r)\bigr)}\int_{B(z,r)} f\, d \lambda
\end{equation*}
 when this integral there exists.  We start from new even for $n=1$ 

\begin{theorem}
  Let $ D\subset \mathbb C^n$ be a pseudoconvex domain, $z_0\in D$, and $u$ be a plurisubharmonic function on $D$, $u (z_0)\neq -\infty$. Then for each number  	$\varepsilon >0$ there is  a holomorphic function $f$ on $D$ such that $f(z_0)\neq 0$ and
$$
	\log \bigl|f(z)\bigr|\leq B_u(z,r)+ n\log \frac{1}{r}+(n+\varepsilon )\log\bigl(1+|z|+r\bigr) 
$$
 for all $z\in D$ and\/ $0<r<\inf \bigl\{|z'-z|\colon z'\in \mathbb C^n\setminus D\bigr\}=:\text{\rm dist\,}(z, \mathbb C \setminus D)$.
  \end{theorem}

This Theorem develops and generalizes results of O.\,V.~Epifanov for $n=1$ [1; Lemma], our results [2; Lemma 1.1], [3, Main Theorem] etc.
The proof uses the H\"ormander--Bombieri method [4; Theorem 4.2.7]. 

\section{{Case $n=1$}}  

The symbol $\mathbb C_{\infty} := \mathbb C \cup  \{\infty\}$  denotes the Riemann sphere, i.\,e. the extended complex plane.
Let $D$ be a subdomain of\/ $\mathbb C_{\infty}\neq D$ everywhere below. We denote by $\text{Hol}(D)$ and $\text{sbh}(D)$ the class of holomorphic  and subharmonic functions on $D$, resp. The class $\text{sbh}(D)$  contains the function $\boldsymbol{-\infty}\colon z\mapsto -\infty$, $z\in D$. 
Let $u\in \text{\rm sbh}(D)\setminus \{\boldsymbol{-\infty}\}$ with the Riesz measure $\nu_u$ hereinafter.

\noindent
\begin{theorem}  Let  $d\colon D\rightarrow (0,+\infty)$ be a function  such that
\begin{equation*}
	0< d(z)< \min \bigl\{\text{\rm dist\,}(z, \mathbb C \setminus D),  1+|z| \bigr\}\quad \text{for all $z\in D$},
\end{equation*}
and  this function is locally separated from zero on $D$, i.\,e. for each  $z\in D$ there is a number $r_z>0$ such that 
	$	B(z,r_z)\subset D$ and $\inf\limits_{z'\in B(z,r_z)} d(z')>0	$. 
If at least one of the following three conditions: 
\begin{description}
	\item[\rm 1)]\label{it:iD}  the closure of  $D$ is not equal to $\mathbb C_{\infty}$;
	\item[\rm 2)]\label{it:ii+D} for the Riesz measure $\nu_u$ of $u$ the inequality  $\nu_u(D)>1$ holds;
	\item[\rm 3)]\label{it:iiD} the domain\/ $D$ is simply connected in\/  $\mathbb C_{\infty}$ and  $\sup\limits_D d<+\infty$ if
	$\nu_u(D)\leq 1$; 
	\end{description}
 is fulfilled,  then there is a non-zero function\/ $f\in \text{\rm Hol}(D)$ such that
\begin{equation*}\label{uvnew}
\log \bigl| f(z)\bigr|\leq B_{u}\bigl (z,d(z)\bigr) +\log\frac1{d(z)}\quad \text{for all\/ $z\in D$}.
\end{equation*}
\end{theorem}

\section{Case $D=\mathbb C$} 

\begin{theorem}  Let  $D=\mathbb C$, $z_0\in D$ and $u(z_0)\neq -\infty$. Then for each number $N \geq 0$ there is  an entire function $f$ such that 
 $$
f(z_0)\neq 0 \quad \text{and}\quad 	\log \bigl|f(z)\bigr|\leq B_u\Bigl(z,  \bigl(1+|z|\bigr)^{-N}\Bigr) \quad \text{for all $z\in \mathbb C$.}
$$
 \end{theorem}

The following  Theorem 4 is a development of H\"ormander's results  [5; {\bf 8}].

\begin{theorem} The following five statements are equivalent:
\begin{description}
	\item[\rm 1)] for every number $N> 0$, there is non-zero entire function  $f$ 	such that
	  	\begin{equation*}
\log \bigl| f(z)\bigr|\leq B_u\left(z,  \bigl(1+|z|\bigr)^{-N}\right) -N\log \bigl(1+|z|\bigr)\quad \text{for all\/ $z\in \mathbb C$};
\end{equation*}
	\item[\rm 2)]  for any number $N> 0$,  there are a positive sequence $r_k\underset{k\to +\infty}{\longrightarrow}+\infty$,  a constant $M>0$,   and a non-zero entire function  $f$ 		such that
	\begin{equation*}
\log \max_{|z|=r_k}\bigl| f(z)\bigr|\leq B_u\bigl(0,  Mr_k\bigr) -N\log \bigl(1+r_k\bigr)\quad \text{for all\/ $k\in \mathbb N$};
\end{equation*}
\item[\rm 3)] $\lim\limits_{r\rightarrow+\infty}\dfrac{1}{\log r}\, B_u(0,r)=+\infty$;
\item[\rm 4)] $\lim\limits_{r\rightarrow+\infty}\dfrac{1}{\log r} \frac{1}{2\pi}\int\limits_{0}^{2\pi} u(re^{i\theta})\, d\theta =+\infty$;
\item[\rm 5)] $\nu_u(\mathbb C)=+\infty$.
\end{description}
\end{theorem}

\section{{Special case $\nu_u(D)<+\infty$}} 
 If the closure of $S\subset D$ is a compact subset of $D\subset \mathbb C_{\infty}$, then  we write $S \Subset D$. For $b\in \mathbb R$, 
we set $b^+:=\max \{0,b\}$. Let\/ $n\in \mathbb N$. Here we assume that  $D\subset \mathbb C_{\infty}$ is  $n$-connected and  $\nu_u(D)<+\infty$.

\begin{theorem}[{[6; Lemma 2.1]}]\label{pr:nn} 
Suppose that $u $ is harmonic, i.\,e.\/ $\nu_u=0$. Then there are a number $b < n-1$ and a function  $f \in  \Hol(D)$ without zeros such that 
	$\log\bigl|f(z)\bigr|\leq  h(z) + b^+ \log \bigl(1+|z|\bigr)$ \quad \text{for all  $z\in D$}.
\end{theorem}

\begin{theorem}
If  	 $0<\nu_u(D)<+\infty$,  then, for each  number  $\varepsilon \in \bigl(0,\nu(D)\bigr)$, there are a constant $b < n-1$, a subdomain   $D_0\Subset D$, and a function  $f \in  \text{\rm Hol}(D)$ without zeros such that 
\begin{multline*}
	\log\bigl|f(z)\bigr|\leq  B_u(z,r) +b^+ \log \bigl(1+|z|\bigr) + \varepsilon \log\frac1{r} 
	\\
	+
	\begin{cases}	\bigl(2\varepsilon -\nu_u(D)\bigr)\log \bigl(1+|z|\bigr) \; &\text{for  $z\in D\setminus D_0$},\\
\nu_u (D) \log \dfrac{1}{r}\;&\text{for  $z\in D_0$},
\end{cases}
\end{multline*}
and for all  $r<\min \bigl\{\text{\rm dist\,}(z, \mathbb C \setminus D),  1+|z| \bigr\}$.
\end{theorem}

\begin{theorem}
If\/  	 $0<\nu_u(D)<+\infty$,  then, for any number  $\varepsilon \in \bigl(0,\nu(D)\bigr)$   and a function $d\colon D\rightarrow (0,+\infty)$ from Theorem\/ {\rm 2}, there are
	a constant $b < n-1$ and a function  $f \in  \text{\rm Hol}(D)$ without zeros such that  
\begin{equation*}
	\log\bigl|f(z)\bigr| \leq  B_u\bigl(z,d(z)\bigr)+b^+ \log\bigl(1+|z|\bigr)  + \varepsilon \log\frac1{d(z)} +
	\bigl(2\varepsilon -\nu_u(D)\bigr)\log \bigl(1+|z|\bigr)  
\end{equation*}
for all $z\in D$. So, for simply  connected domain $D$ we have $b^+=0$.
\end{theorem}

\bigskip
\centerline{\large \bf References}
\medskip
 
\begin{enumerate}

\item {\sl О.\,В.~Епифанов,\/} {О разрешимости неоднородного уравнения Коши--Римана в классах функций, ограниченных с весом и системой весов~/\!/ 
 Матем. за\-м\-е\-т\-ки.---1992.---Т.~51, Вып.~1.---С.~83--92. --- {\sl O. V. Epifanov,\/} On solvability of the nonhomogeneous Cauchy--Riemann equation in classes of functions that are bo\-u\-n\-d\-ed with weights or systems of weights, Mathematical Notes, 1992, 51:1, 54--60.}

\item {\sl Б.\,Н.~Хабибуллин,  Ф.\,Б.~Хабибуллин,  Л.\,Ю.~Чередникова,\/} {Под\-по\-с\-л\-е\-д\-о\-в\-а\-т\-е\-ль\-ности нулей для классов голоморфных функций, их устойчивость и энтропия линейной связности. I~/\!/        Алгебра и анализ. ---2008.---Т.~20, №~1.---С.~146--189. ---
{\sl B. N. Khabibullin, F. B. Khabibullin, L. Yu. Cherednikova,\/} Zero subsequences for classes of holomorphic functions: stability and the entropy of arcwise connectedness. I, St. Petersburg Math. Journal, 2009, 20:1, 101--129.}

\item {\sl  Б.\,Н.~Хабибуллин,\/} {Множества единственности в пространствах целых функций одной переменной~/\!/  Изв. АН СССР. Сер. 
матем.---1991.---Т.~55, №~5.---С.~1101--1123. --- {\sl B. N. Khabibullin,\/} Sets of uniqueness in spaces of entire functions of a single variable, Mathematics of the USSR-Izvestiya, 1992, 39:2, 1063--1084.}

 \item {\sl Lars H\"ormander,\/} {Notions of Convexity---Boston: Birkh\"aser, 1994.}

\item {\sl Lars H\"ormander,\/} {On the Legendre and Laplace transformations~/\!/  Ann. Sc. Norm. Super. Pisa, Cl. Sci., IV. Ser.---1997.---V.~XXV, No.~3--4.---P.~517--568.}

\item {\sl Б.\,Н.~Хабибуллин,\/}  {Последовательности нулей голоморфных функций, представление мероморфных функций и гармонические 
миноранты~/\!/ Матем. сб.---2007.---Т.~198, №~2.---С.~121--160. --- {\sl B.\,N.~Khabibullin,\/} Zero sequences of holomorphic functions, representation of meromorphic functions, and harmonic minorants, Sbornik: Mathematics, 2007, 198:2, 261--298.}
 \end{enumerate}
\end{document}